\documentclass[12pt]{article}
\usepackage{amsfonts}
\usepackage{latexsym,amsmath,amssymb}
\usepackage{graphics,epstopdf}
\usepackage[pdftex]{graphicx}

\numberwithin{equation}{section}
\newtheorem{thm}{Theorem}[section]

\newtheorem{rem}[thm]{Remark}
\numberwithin{equation}{section}
\begin{document}

\bigskip

\bigskip

\begin{center}

{\Large \textbf{Some Approximation Results by $(p, q)$-analogue of
		Bernstein-Stancu Operators (Revised)}}
\bigskip\\

\textbf{M. Mursaleen}, \textbf{Khursheed J. Ansari} and \textbf{Asif Khan}

Department of\ Mathematics, Aligarh Muslim University, Aligarh--202002, India%
\\[0pt]

mursaleenm@gmail.com; ansari.jkhursheed@gmail.com; asifjnu07@gmail.com \\[0pt%
]

\bigskip

\bigskip

\textbf{Abstract}
\end{center}

\parindent=8mm {\footnotesize {In this paper,   We have given a corrigendum to our paper ``Some Approximation Results by $(p,q)$-analogue of Bernstein-Stancu Operators" published in Applied Mathematics and Computation 264 (2015) 392-402. We introduce a new analogue of Bernstein-Stancu operators and we call it as $(p, q)$-Bernstein-Stancu operators. We study approximation properties based on Korovkin’s type approximation theorem of $(p, q)$-Bernstein-Stancu operators. We also establish some direct theorems. }}\newline



\section{Construction of Revised Operators}

Mursaleen et. al \cite{mka3} introduced $(p,q)$-analogue of Bernstein-Stancu operators as
\begin{equation*}
S_{n,p,q}(f;x)=\sum\limits_{k=0}^{n}\left[
\begin{array}{c}
n \\
k%
\end{array}%
\right] _{p,q}x^{k}\prod\limits_{s=0}^{n-k-1}(p^{s}-q^{s}x)~~f\left( \frac{%
[k]_{p,q}+\alpha}{[n]_{p,q}+\beta}\right) ,~~x\in \lbrack 0,1].\eqno(1)
\end{equation*}
But $S_{n,p,q}(1;x)\ne1$ for all $x\in[0,1)$ and $S_{n,p,q}(1;x)=1$ for $x=1$ only. Hence, we are re-introducing our operators as follows:
\begin{equation*}
S_{n,p,q}(f;x)=\frac1{p^{\frac{n(n-1)}2}}\sum\limits_{k=0}^{n}\left[
\begin{array}{c}
n \\
k%
\end{array}%
\right] _{p,q}p^{\frac{k(k-1)}2}x^{k}\prod\limits_{s=0}^{n-k-1}(p^{s}-q^{s}x)~~f\left( \frac{%
{p^{n-k}[k]_{p,q}+\alpha}}{[n]_{p,q}+\beta}\right) ,~~x\in \lbrack 0,1].\eqno(2)
\end{equation*}

Note that for $\alpha=\beta=0,$ $(p,q)$-Bernstein-Stancu operators given by (2) reduces into $(p,q)$-Bernstein operators as given in \cite{mka1,mka2}.
 Also for $p=1$, $(p,q)$-Bernstein-Stancu operators given by (2) turn out to
be $q$-Bernstein-Stancu operators.\\

For any $p>0$ and $q>0,$ the $(p,q)$ integers $[n]_{p,q}$ are defined by

\begin{equation*}
\lbrack n]_{p,q}=p^{n-1}+p^{n-2}q+p^{n-3}q^2+...+pq^{n-2}+q^{n-1}\\
=\left\{
\begin{array}{lll}
\frac{p^{n}-q^{n}}{p-q},~~~~~~~~~~~~~~~~\mbox{when $~~p\neq q \neq 1$  } & \\
&  \\
n~p^{n-1},~~~~~~~~~~~~~~\mbox{ when $p=q\neq1$  } & \\
&  \\

[n]_q ,~~~~~~~~~~~~~~~~~~~\mbox{when $p=1$  }& \\

n ,~~~~~~~~~~~~~~~~~~~~~\mbox{ when $p=q=1$  }
\end{array}%
\right.
\end{equation*}
~where  $[n]_q $ denotes the $q$-integers and $n=0,1,2,\cdots$.\\

Obviously, it may be seen that $[n]_{p,q}= p^{n-1}[n]_{\frac{q}{p}}.\\$

We have the following basic result:\newline

\parindent=0mm\textbf{Lemma 1.} For $x\in \lbrack 0,1],~0<q<p\leq 1$, we
have\newline
\newline
(i)~~$S_{n,p,q}(1;x)=~1$,\newline
(ii)~$S_{n,p,q}(t;x)=~\frac{[n]_{p,q}}{[n]_{p,q}+\beta}x+\frac{\alpha}{[n]_{p,q}+\beta}$,\newline
(iii) $S_{n,p,q}(t^{2};x)=\frac{q[n]_{p,q}[n-1]_{p,q}}{([n]_{p,q}+\beta)^2}x^2+\frac{[n]_{p,q}(2\alpha+p^{n-1})}{([n]_{p,q}+\beta)^2}x+\frac{\alpha^2}{([n]_{p,q}+\beta)^2}$.\newline

\parindent=0mm\textbf{Proof}. (i)
\begin{equation*}
S_{n,p,q}(1;x)=\frac1{p^{\frac{n(n-1)}2}}\sum\limits_{k=0}^{n}\left[
\begin{array}{c}
n \\
k%
\end{array}%
\right] _{p,q}p^{\frac{k(k-1)}2}x^{k}\prod\limits_{s=0}^{n-k-1}(p^{s}-q^{s}x)=1.
\end{equation*}

(ii)
\begin{eqnarray*}
S_{n,p,q}(t;x) &=&\frac1{p^{\frac{n(n-1)}2}}\sum\limits_{k=0}^{n}\left[
\begin{array}{c}
n \\
k%
\end{array}%
\right] _{p,q}p^{\frac{k(k-1)}2}x^{k}\prod\limits_{s=0}^{n-k-1}(p^{s}-q^{s}x)~~\frac{%
p^{n-k}[k]_{p,q}+\alpha}{[n]_{p,q}+\beta}\\
&=&\frac1{p^{\frac{n(n-1)}2}}\sum\limits_{k=0}^{n}\left[
\begin{array}{c}
n \\
k%
\end{array}%
\right] _{p,q}p^{\frac{k(k-1)}2}x^{k}\prod\limits_{s=0}^{n-k-1}(p^{s}-q^{s}x)~~\frac{%
p^{n-k}[k]_{p,q}}{[n]_{p,q}+\beta}\\
~~&+&\frac{\alpha}{[n]_{p,q}+\beta}\frac1{p^{\frac{n(n-1)}2}}\sum\limits_{k=0}^{n}\left[
\begin{array}{c}
n \\
k%
\end{array}%
\right] _{p,q}p^{\frac{k(k-1)}2}x^{k}\prod\limits_{s=0}^{n-k-1}(p^{s}-q^{s}x)\\
&=&\frac{[n]_{p,q}}{p^{\frac{n(n-3)}2}}\sum\limits_{k=0}^{n-1}\left[
\begin{array}{c}
n-1 \\
k%
\end{array}%
\right] _{p,q}p^{\frac{(k+1)(k)}2}x^{k+1}\prod\limits_{s=0}^{n-k-2}(p^{s}-q^{s}x)\frac{1}{p^{k+1}([n]_{p,q}+\beta)}+\frac{\alpha}{[n]_{p,q}+\beta}\\
&=&\frac {[n]_{p,q}}{[n]_{p,q}+\beta}x ~~\frac{1}{p^{\frac{(n-1)(n-2)}{2}}} \sum \limits_{k=0}^{n-1}\left[
\begin{array}{c}
n-1 \\
k%
\end{array}%
\right] _{p,q}p^{\frac{k(k-1)}2}x^{k}\prod\limits_{s=0}^{n-k-2}(p^{s}-q^{s}x)+\frac{\alpha}{[n]_{p,q}+\beta}\\
&=&\frac{[n]_{p,q}}{[n]_{p,q}+\beta}x+\frac{\alpha}{[n]_{p,q}+\beta}.
\end{eqnarray*}

(iii)
\begin{eqnarray*}
S_{n,p,q}(t^2;x) &=&\frac1{p^{\frac{n(n-1)}2}}\sum\limits_{k=0}^{n}\left[
\begin{array}{c}
n \\
k%
\end{array}%
\right] _{p,q}p^{\frac{k(k-1)}2}x^{k}\prod\limits_{s=0}^{n-k-1}(p^{s}-q^{s}x)~~{\bigg(\frac{p^{n-k}[k]_{p,q}+\alpha}{[n]_{p,q}+\beta}\bigg)}^2\\
&=&\frac{1}{([n]_{p,q}+\beta)^2}~\frac1{p^{\frac{n(n-1)}2}}\Bigg[p^{2n}\sum\limits_{k=0}^{n}\left[
\begin{array}{c}
n \\
k%
\end{array}%
\right] _{p,q}p^{\frac{k(k-1)}2}x^{k}\prod\limits_{s=0}^{n-k-1}(p^{s}-q^{s}x)~\frac{[k]_{p,q}^2}{p^{2k}}\\
&&+2\alpha ~p^n\sum\limits_{k=0}^{n}\left[
\begin{array}{c}
n \\
k%
\end{array}%
\right] _{p,q}p^{\frac{k(k-1)}2}x^{k}\prod\limits_{s=0}^{n-k-1}(p^{s}-q^{s}x)~\frac{[k]_{p,q}}{p^{k}}\\
&&+\alpha^2~\sum\limits_{k=0}^{n}\left[
\begin{array}{c}
n \\
k%
\end{array}%
\right] _{p,q}p^{\frac{k(k-1)}2}x^{k}\prod\limits_{s=0}^{n-k-1}(p^{s}-q^{s}x)\Bigg].\\
\end{eqnarray*}
On shifting the limits and using $[k+1]_{p,q}=p^k+q[k]_{p,q}$, we get our desired result.
\newline

\parindent=0mm\textbf{Lemma 2.2}. For $x\in [0,1],~0<q<p\leq 1$\\
\\
(i)~~~$S_{n,p,q}\bigl{(}(t-x);x\bigl{)}=\frac{\alpha-\beta x}{[n]_{p,q}+\beta}$,\\
(ii)~~$S_{n,p,q}\bigl{(}(t-x)^2;x\bigl{)}= \big\{ \frac{q[n]_{p,q}[n-1]_{p,q}-[n]_{p,q}^2+{\beta}^2}{([n]_{p,q}+\beta)^2}\big\}x^2+\big\{ \frac{p^{n-1}[n]_{p,q}-2 \alpha \beta}{([n]_{p,q}+\beta)^2} \big\}x+\frac{\alpha^2}{([n]_{p,q}+\beta)^2}$.\\
\\
\parindent=0mm\textbf{Proof}. (ii) By linearity of the operator we have
\begin{eqnarray*}
 S_{n,p,q}\bigl{(}(t-x)^2;x\bigl{)} &=&S_{n,p,q}(t^2;x)-2xS_{n,p,q}(t;x)+x^2S_{n,p,q}(1;x) \\
   &=&\frac{q[n]_{p,q}[n-1]_{p,q}}{([n]_{p,q}+\beta)^2}x^2+\frac{[n]_{p,q}(2\alpha+p^{n-1})}{([n]_{p,q}+\beta)^2}x\\
   &&+\frac{\alpha^2}{([n]_{p,q}+\beta)^2}-2x \frac{[n]_{p,q}}{[n]_{p,q}+\beta}-2x\frac{\alpha}{[n]_{p,q}+\beta}+x^2\\
&=& \bigg\{ \frac{q[n]_{p,q}[n-1]_{p,q}-[n]_{p,q}^2+{\beta}^2}{([n]_{p,q}+\beta)^2}\bigg\}x^2+\bigg\{ \frac{p^{n-1}[n]_{p,q}-2 \alpha \beta}{([n]_{p,q}+\beta)^2} \bigg\}x+\frac{\alpha^2}{([n]_{p,q}+\beta)^2}.
\end{eqnarray*}

\section{ Main Results:}

\subsection{Korovkin type approximation theorem}

\begin{rem}\label{r5.1}
	For $q\in(0,1)$ and $p\in(q,1]$, it is obvious that
	$\lim\limits_{n\to\infty}[n]_{p,q}=0 $ or $\frac1{p-q}$. In order to reach to convergence
	results of the operator $L^{n}_{p,q}(f;x),$ we take a sequence $q_n\in(0,1)$ and $p_n\in(q_n,1]$
	such that $\lim\limits_{n\to\infty}p_n=1,$ $\lim\limits_{n\to\infty}q_n=1$ and $\lim\limits_{n\to\infty}p_n^n=1,$ $\lim\limits_{n\to\infty}q_n^n=1$. So we get
	$\lim\limits_{n\to\infty}[n]_{p_n,q_n}=\infty$.
\end{rem}

%
%

\parindent=0mm\textbf{Theorem 3.1.1}. Let sequence $q_n$ and $p_n $ be such that it satisfies Remarks $\ref{r5.1}$, then $S_{n,p,q}(f,x)$   converges uniformly to $f$ on $[0,1]$.\\
\\
\parindent=0mm\textbf{Proof}. By the Bohman-Korovkin Theorem it is sufficient to show that
$$\lim\limits_{n\to\infty}\|S_{n,p_n,q_n}(t^m;x)-x^m\|_{C[0,1]}=0,~~~m=0,1,2.$$
By Lemma 2.1 (i)-(ii), it is clear that
\begin{equation*}
\lim\limits_{n\to\infty}\|S_{n,p_n,q_n}(1;x)-1\|_{C[0,1]}=0;
\end{equation*}
\begin{equation*}
\lim\limits_{n\to\infty}\|S_{n,p_n,q_n}(t;x)-x\|_{C[0,1]}=0.
\end{equation*}
Using $q_n[n-1]_{p_n,q_n}=[n]_{p_n,q_n}-p_n^{n-1}$ and by Lemma 2.1(iii), we have
\begin{equation*}
  |S_{n,p_n,q_n}(t^2;x)-x^2|_{C[0,1]}=\biggl{|}\frac{q[n]_{p,q}[n-1]_{p,q}}{([n]_{p,q}+\beta)^2}x^2
  +\frac{[n]_{p,q}(2\alpha+p^{n-1})}{([n]_{p,q}+\beta)^2}x+\frac{\alpha^2}{([n]_{p,q}+\beta)^2}-x^2\biggl{|}.
\end{equation*}

Taking maximum of both sides of the above inequality, we get
\begin{equation*}
\|S_{n,p_n,q_n}(t^2;x)-x^2\|_{C[0,1]}\leq\frac{q[n]_{p,q}[n-1]_{p,q}}{([n]_{p,q}+\beta)^2}-1
+\frac{[n]_{p,q}(2\alpha+p^{n-1})}{([n]_{p,q}+\beta)^2}+\frac{\alpha^2}{([n]_{p,q}+\beta)^2}
\end{equation*}
which yields
\begin{equation*}
\lim\limits_{n\to\infty}\|S_{n,p_n,q_n}(t^2;x)-x^2\|_{C[0,1]}=0.
\end{equation*}
\\
Thus the proof is completed and all other results follow similarly.\\



\begin{thebibliography}{99}
\bibitem{mka1} M. Mursaleen, Khursheed J. Ansari, Asif Khan, On $(p,q)$-analogue of Bernstein operators, Applied Mathematics and Computation, 266 (2015) 874-882.

    \bibitem{mka2} M. Mursaleen, Khursheed J. Ansari, Asif Khan, On $(p,q)$-analogue of Bernstein operators (revised), \emph{arxiv}:1503.07404v2[math CA].

\bibitem{mka3} M. Mursaleen, K. J. Ansari and Asif Khan,, Some Approximation Results by $(p,q)$-analogue of Bernstein-Stancu Operators, \emph{Applied Mathematics and Computation} 264,(2015), 392-402.


\end{thebibliography}
\end{document}